\documentclass[letterpaper, 10 pt, conference]{ieeeconf}  

\IEEEoverridecommandlockouts                              

\usepackage{color}
\usepackage{graphics} 
\usepackage{epsfig} 
\usepackage{subfigure}
\usepackage{times} 
\usepackage{amsmath} 
\usepackage{amssymb}  
\usepackage{lipsum}
\usepackage{multirow}

\newtheorem{definition}{Definition}[section]

\newtheorem{remarkth}[definition]{Remark}

\usepackage{amssymb}

\newcommand{\proa}{A^*G \mbox{$\;$}_{\tau^*} \kern-3pt\times_\alpha
G \mbox{$\;$}_\beta \kern-3pt\times_{\tau^*} A^*G}

\title{\LARGE \bf
Optimal Control of Quantum Purity for $n=2$ Systems
}

\author{William Clark$^{1}$,  Anthony Bloch$^{1}$, Leonardo Colombo$^{1}$ and Patrick Rooney$^{2}$
\thanks{$^{1}$W. Clark, A. Bloch, and L. Colombo are with Department of Mathematics, University of Michigan, 530 Church St. Ann Arbor, 48109, Michigan, USA.
        {\tt\small wiclark@umich.edu, abloch@umich.edu, ljcolomb@umich.edu}}%
\thanks{$^{2}$ P. Rooney is with Department of Physics, University of Windsor, Ontario N9B 3P4, Canada.
	    {\tt\small darraghrooney@gmail.com}}
}

\begin{document}

\maketitle
\thispagestyle{empty}
\pagestyle{empty}

\begin{abstract}
The objective of this work is to study time-minimum and energy-minimum global optimal control for dissipative open quantum systems whose dynamics is governed by the Lindblad equation. The controls appear only in the Hamiltonian.

Using recent results regarding the decoupling of such dissipative dynamics into intra- and inter-unitary orbits, we transform the control system into a bi-linear control system on the Bloch ball (the unitary sphere together with its interior). We then  design a numerical algorithm to construct an optimal path to achieve a desired point given initial states close to the origin (the singular point) of the Bloch ball. This is done both for the minimum-time and minimum -energy control problems. 

\end{abstract}

\section{Introduction}
Control of quantum conservative (Hamiltonian) systems has been extensively studied in the last few decades from both theoretical and interdisciplinary points of view \cite{BCS}, \cite{2-3level}, \cite{Brocket1}, \cite{Brocket2}. Recently, there has been a growing interest in control of open (dissipative, non-Hamiltonian) quantum
systems  because of their  applications to physics, chemistry
and quantum computing. For example there has been interest in the control of the rotation of a molecule in the gaseous phase by using laser fields in dissipative media \cite{RS}. Here the dissipation is due to molecular collisions and the dimension of the Hilbert space describing the states of the system (infinite-dimensional) can be truncated to make  the system finite-dimensional if the intensity of the laser field is sufficiently weak. Other applications include control of the spin dynamics by magnetic fields in nuclear magnetic resonance \cite{EB} and applications to the  construction of quantum computers \cite{RB}. 

The aim of this work is to study optimal control of two-level quantum systems in a  dissipative environment, where we assume that the dissipation is Markovian (the dynamics depends only on the present state and not its history) and time-independent.  In this case the evolution for the density matrix of the system can by described by a quantum dynamical semi-group and the Lindblad master equation \cite{Altafini}, \cite{Breuer}, \cite{lindblad1976}. 

The state space for a \textit{closed quantum system} is an $n$-dimensional projective Hilbert space, $P(\mathcal{H})$, of a complex Hilbert space $\mathcal{H}$. Typically one drops the requirement of the space being projective, and instead, we work with unit vectors. To preserve the length of the vectors, \textit{time evolution is unitary} $\displaystyle{U(t_1,t_2)|\psi(t_1)\rangle = |\psi(t_2)\rangle}$, and the evolution is described by the \textit{Schr\"{o}dinger equation}
\begin{equation}
\frac{d}{dt}|\psi(t)\rangle = -iH(t)|\psi(t)\rangle,
\end{equation}
where $H$ is the Hermitian Hamiltonian. Here the \textit{ket}-bracket describes the vector associated with an observable state.






The \textit{density operator}, $\rho$, describes a probabilistic ensemble of states. It is given by a positive semi-definite Hermitian operator $\rho$ with $\text{Tr}(\rho)=1$ and $\text{Tr}(\rho^2)\leq 1$. The \textit{purity} of a density operator describes how close $\rho$ is to a single state. It is usually defined as
$\displaystyle{P_2(\rho)={\text{Tr}(\rho^2)}\in [1/n,1]}$, where the unique operator that has a purity of $1/n$ is $\frac{1}{n}I_{n\times n}$, called the \textit{completely mixed state.}  The dynamics for the purity operator is described by the \textit{von Neumann equation}
\begin{equation}
\frac{d}{dt}\rho = [-iH,\rho].
\end{equation}
Notice that this dynamical system preserves the purity of $\rho$ (because the system is iso-spectral). A consequence is that if the quantum system is controlled by its Hamiltonian, there is no controllability over its purity, since one cannot directly alter the probabilities or achieve a purity of one.


\textit{Open quantum systems} are quite different. For such systems dissipation occurs when we allow the system to interact with the environment. The full picture is an integro-differential equation called the Nakajima-Zwanzig (NZ) equation. To make the dissipation purely of differential form, one usually make two assumptions: the dissipation is Markovian (i.e. the dissipation only depends upon the current state, not past history) and the dissipation is time-invariant.

Under these assumptions, the dynamics of the density operator is given by the \textit{Lindblad master equation} \cite{lindblad1976}
$$\frac{d}{dt}\rho = [-iH,\rho]+\sum_{j=1}^{N}\left(
L_j\rho L_j^{\dagger} - \frac{1}{2}\left\{ L_j^{\dagger}L_j,\rho\right\}\right)
$$ where the $L_j$ are called the Lindblad operators, $N$ denotes the quantity of Lindblad operators, and $\{\cdot,\cdot\}$ is the anti-commutator: $\{A,B\}=AB+BA$. Along the work we assume Linblad operators are traceless.

\subsection{Main contributions:} In this paper we investigate the problem of minimum-time and minimum-energy global optimal control for dissipative open quantum systems whose dynamics is governed by the Lindblad equation. Our contributions here are two-fold. First we improve and extend the results for local optimality based on the steepest descent method studied in \cite{PhysRevA.93.063424} by obtaining global results on the Bloch ball, which is the physical state space of the system. Also we consider an  energy-minimum optimal control problem, where the cost corresponds to the energy transfer between the control and the internal Hamiltonian. The second contribution is related to work on generalizing  the results given in  \cite{BonnardJMP}, \cite{Bonnarcyber}, \cite{bonnarsiam}, \cite{Sugny} with bounded controls to a class of control system with more general Lindblad operators. 

\subsection{Outline:} The structure of the work is as follows: In Section II we introduce the Lindblad equation and we interpret its dynamics as a control system in the Bloch ball. Section III explains why the optimal control problem is singular and how to achieve maximum purity by suitable choice of initial conditions for the boundary value problem. Sections IV and V are devoted to the study of time-minimum and energy-minimum controls for two- and three-dimensional systems. Numerical results for time-minimum and energy-minimum controls in the previous two situations are explored in Section VI. We conclude in Section VII by outlining future research.
\section{The Lindblad equation}
\label{section2} 

An open quantum system is described by a density operator $\rho$, which is a trace-one positive semi-definite Hermitian operator on an $n$-dimensional complex Hilbert space $\mathcal{H}$. If the dissipation is Markovian and time independent, the density operator obeys the Lindblad equation
(see \cite{Breuer} and \cite{lindblad1976} for details) 
\begin{equation}\label{eq:master}
\frac{d\rho}{dt} = [-iH,\rho]+\mathcal{L}_D(\rho),
\end{equation} with \begin{equation}\label{eq:Lindblad}
\mathcal{L}_D(\rho) := \sum_{j=1}^{N} L_j\rho L_j^{\dagger} - \frac{1}{2}\{L_j^{\dagger}L_j,\rho\}.
\end{equation} where $N$ denotes the quantity of Lindblad operators,
$[\cdot,\cdot]$ denotes commutator of matrices, $H$ is the Hermitian Hamiltonian, $\dagger$ represents the Hermitian transpose and $\left\{L_j\right\}$ are the Lindblad operators. The purity of the system is defined as $P_{2}(\rho)={\text{Tr}(\rho^2)}$.


The goal is to construct controls for the purity operator under Lindblad dissipation. In most situations the controls appear in the Hamiltonian $H$, and not in the Lindblad operator  $\mathcal{L}_{D}$. This is the assumption we make here.

\indent When $n=2$, the density operator can be identified with a vector in the Bloch ball (the unitary sphere with its interior) \cite{5288564}. Under this special case, we can change the view from dynamics on operators to dynamics in the Bloch ball by considering $$\rho=\frac{1}{2}\left(I+\sum_{j=1}^{3}q_j\sigma_j\right)$$ where $q\in S^{2}$ (i.e., $q_1^2+q_2^2+q_3^2\leq 1$), $I$ is the $2\times 2$ identity matrix, and $\sigma_{j}$ are the Pauli matrices.

Using this identification, we can reformulate the Lindblad equation (\ref{eq:master}) into a first-order dynamical system on the unit ball. Using the derivation given in \cite{PhysRevA.93.063424} (see appendix A therein), we have that \eqref{eq:master} is equivalent to
\begin{equation}\label{eq:bloch}
\frac{d\vec{q}}{dt} = \vec{b} + (A-\text{tr}(A))\vec{q} + \vec{u}\times\vec{q},
\end{equation}
where 
\begin{equation}\label{eq:Amatrix}
A :=\frac{1}{2}\sum_j \vec{l_j}\overline{\vec{l}_j}^T+\overline{\vec{l}_j}\vec{l}_j^T,\qquad
\vec{b} := i\sum_j\! \vec{l}_j \times \overline{\vec{l}_j},
\end{equation}
given that the bar represents the complex conjugate of matrices.
The vectors $l_j$, $u$ are the traceless parts of $L_j$ and $H$ respectively with $H = \displaystyle{h_0I+\sum_{k=1}^3 \! u_k\sigma_k}$ where $\sigma_k$ are the Pauli matrices. Notice that the matrix $A$ is positive semi-definite. 
From here on out, we will call the matrix $A-\text{tr}(A)$ to be $B$.

Here we would like to point out that the system \eqref{eq:bloch} reduces in special cases to those  studied in \cite{BonnardJMP}, \cite{Bonnarcyber}, \cite{bonnarsiam}, \cite{Sugny}.  Setting one of the controls in our system to zero and  identifying the parameters that appear in the system given in  \cite{Bonnarcyber}, \cite{Sugny}  with the elements of the matrix A,
one obtains the system discussed in those papers.

As discussed above the control variables in the open 
quantum systems we discuss here appear in the Hamiltonian operator as
in \cite{5288564}, \cite{Sugny}. The controlled Hamiltonian dynamics cannot achieve a purity one \cite{tannor} and in general  cannot affect the purity of the state or transfer the states between unitary orbits. To control purity one must use the dissipative dynamics to move between orbits as in \cite{PhysRevA.93.063424} and \cite{Rooney2}.  In this paper  we consider \textit{unbounded} controls $\left\{u_{k}\right\}$ which may take any value in $\mathbb{R}$.

As long as $q$ is not at the origin, it can be shown that equation \eqref{eq:bloch} can be written in terms of the radial component (i.e., $r=\lVert q\rVert$). In \cite{PhysRevA.93.063424}  it is shown that the purity $P_{2}(\rho)$ is equal to $\displaystyle{\frac{(1+r^2)}{2}}$, where $r=\lVert q\rVert$. So, controlling the purity is synonymous with controlling the magnitude of the Bloch vector. It is  be helpful to extract the dynamics for $r$.  Knowing that $r^2=\langle q,q\rangle$, we find that
\begin{equation*}
\begin{split}
2r\dot{r} &= 2\langle \dot{q},q\rangle
= 2\langle b+u\times q + Bq,q\rangle = 2\langle b,q\rangle +2\langle Bq,q\rangle\\
&= 2r\langle b,\hat{q}\rangle + 2r^2\langle B\hat{q},\hat{q}\rangle
\end{split}
\end{equation*} where $\hat{q}$ is the unit vector associated with $q$, $\hat{q} = q /||q||$.

Therefore, 
\begin{equation}\label{eq:redial}
r\frac{dr}{dt} = \langle q, b+Bq\rangle := f(q).
\end{equation} So, we can control the purity by controlling the orientation of the corresponding unit vector.  Hereafter we refer to $\displaystyle{f(q)}$ as \textit{the purity derivative}.

Our goal is to find a control scheme that optimally transports the completely mixed state to a state of maximal purity. This raises two questions: What is the maximal achievable purity? and, what do we mean by optimal?. To answer these questions we will introduce in the next section the notion of \textit{apogee} and \textit{escape chimney} as in \cite{PhysRevA.93.063424}.


\section{The apogee and the escape chimney} \label{section3}
The purity derivative (\ref{eq:redial}) is independent of the controls used. It is illuminating to examine the regions in the Bloch sphere where the purity derivative is positive. To do this, we examine the zeros of $f$. Define two sets, $\mathcal{U} = \{q|f(q)\geq 0\}$ and the ellipsoid $\mathcal{M} = \{q|f(q)=0\}$. To find $\mathcal{M}$, we define a new function $f_q(r) := f(rq)$ where $q\in S^2$. Finding the roots of $f_q$ will let us solve for $\mathcal{M}$ in spherical coordinates.
\begin{equation}\label{eq:radialzeros}
f_q(r) = \langle q,Bq\rangle r^2 + \langle q,b\rangle r.
\end{equation}
So, the nonzero root is
\begin{equation}\label{eq:nonzeroroot}
g(q) := -\frac{\langle q,b\rangle}{\langle q,Bq\rangle}.
\end{equation}
Notice that $g$ is always defined since $B$ is negative-definite, and also note that the maximum of $g$ must be bounded by $1$, so the Bloch ball is invariant under (\ref{eq:bloch}).\\
\indent Define $q_{apogee}$ to be \textit{the apogee} of the ellipsoid $\mathcal{M}$. i.e.
\begin{equation}\label{eq:endpoint}
q_{apogee} := \arg\max_{q\in\mathcal{M}} ~ \lVert q \rVert,
\end{equation}
which can be found by maximizing \eqref{eq:nonzeroroot} on $S^2$. Therefore, the apogee of the ellipsoid $\mathcal{M}$ will be the state with the maximal achievable purity. We call the interior of the ellipsoid $\mathcal{U}$ the \textit{escape chimney}. 

In Figure 3, in Section VI-B, we show a picture of the escape chimney inside the Bloch ball where the black square is the apogee and the trajectories represents energy minimum and time minimum control for initial conditions near the origin of the Bloch ball and final conditions at the apogee.

 The optimal control problem consists of finding a trajectory of the state variables, starting at the completely mixed state (i.e., $r=0$) and ending  at the apogee. 
It is important to note, however, that the dynamics (\ref{eq:redial})  has a singularity at the origin since $$\frac{dr}{dt}=\frac{f(q)}{r},$$ as well as the fact that the apogee cannot be reached in finite time since it is not possible to reach equilibrium points in finite time. Note that the purity derivative is independent of controls, so that, for a given path the purity derivative is an autonomous first order dynamical system that cannot reach its fixed point.

To circumvent these problems, we take the following boundary conditions representing initial and final states
\begin{equation}\label{eq:startandend}
q_0 = \varepsilon \vec{b}/\lVert b\rVert, \hspace{0.15in} q_f = (1-\delta)q_{apogee}
\end{equation}
with $\varepsilon>0$ and $\delta>0$ sufficiently small.

Therefore, we can now state the optimal control problem studied in this work as follow: Let $J:S^2\times U\rightarrow\mathbb{R}$ be a cost functional dependent on the state as well as the controls. The optimal control problem consists of finding a control, $u:[0,t_f]\rightarrow U$ satisfying the dynamics (\ref{eq:redial}) such that $q(0)=q_0$, $q(t_f)=q_f$ and $u = \arg\min\! \int_0^{t_f} \! J(q,u)\, dt.$ In this work we study two different optimal control problems depending on the cost functional we choose: A \textit{time-minimal optimal control problem} ($J=1$) and an \textit{energy-minimal optimal control problem} ($J=||u||_{2}^{2}$).


\section{Two-dimensional Systems}
In the special case when $N=1$ in (\ref{eq:Lindblad}), that is, only one Lindblad term; $\vec{b}$ becomes an eigenvector of $B$. This fact lets us simultaneously diagonalize $B$ and rotate $\vec{b}$ into the first coordinate. By additionally taking $u_1=u_2=0$, the third component of $q$ in equation \eqref{eq:bloch} becomes uncontrolled and exponentially decays to zero. Dropping this coordinate, our system collapses to a two-dimensional underactuated bi-linear control system:
\begin{equation}\label{eq:twodimenions}
\begin{array}{rcl}
\dot{x} &=& b_1 + a_1 x -uy \\
\dot{y} &=& b_2 + a_2 y + ux
\end{array} 
\end{equation} where $u=u_3$, $a_1,a_2$ are the coefficients of the matrix $A$, $q=(x,y)$ and $\vec{b}=(b_1,b_2)$.
\subsection{Time-Minimal Controls}
We want to find (unbounded) controls that steer (\ref{eq:twodimenions}) with end points (\ref{eq:startandend}) in the minimal amount of time possible. i.e. find a minimal solution to the functional
\begin{equation}\label{eq:2dlagrangian}
\min \int_{0}^{t_f} \! dt = \min \int_{x_0}^{x_f} \! \frac{dt}{dr}\frac{dr}{dx} \, dx
\end{equation} where $x(0)=x_0$ and $x(t_f)=x_f$.
To find $dr/dt$, see $(\ref{eq:redial})$ and $dr/dx = x+yy'$. So we wish to minimize a functional with integrand 
\begin{equation}\label{eq:lagrangian}
I(q,q') = \int\! L(q,q')  dt = \int_{x_0}^{x_f} \! \frac{x+yy'}{\langle q,b+Bq\rangle} \, dx.
\end{equation}
This Lagrangian is not hyperregular, so the Euler-Lagrange equations will fail to yield meaningful results \cite{baillieul2008nonholonomic}.
To get around this problem, we implement the Rayleigh-Ritz numerical algorithm \cite{hoffman2001numerical}. We assume $y(x)$ is a sum of linearly independent functions
\begin{equation}\label{eq:yform}
y(x) = y^0(x) + \sum_{i=1}^M c_iy^i(x),
\end{equation}
where $y^0(x_0) = y_0$, $y^0(x_f)=y_f$, and $y^i(x_0)=y^i(x_f)=0$. Specifically, we will take the following functions as a basis of polynomials for our approximation.
\begin{equation}\label{eq:guessfunctions}
\begin{split}
&y^i(x) = (x-x_0)(x-x_f)^i,\quad i=1,\ldots,M;\\
&y^0(x) = \frac{y_f-y_0}{x_f-x_0}(x-x_0)+y_0,
\end{split}
\end{equation} with $M$ an arbitrary integer.
Then, a necessary condition for our guess to minimize the functional \eqref{eq:lagrangian} is for the following $M$ equations to hold
\begin{equation}\label{eq:partials}
\frac{\partial}{\partial c_i} I = \int_{x_0}^{x_f}\! \frac{\partial}{\partial c_i} L \, dx = 0.
\end{equation}
This can be done by symbolically computing $\partial L/ \partial c_i$ in MATLAB and numerically integrating using a $4^{th}$ order Runge-Kutta method. In order to find the optimal values to the $c_i$'s, we construct a new function
\begin{eqnarray}\label{eq:norm}
\nu:\mathbb{R}^M&\longrightarrow&\mathbb{R}\nonumber\\
c&\longmapsto &\left( \sum_{i=1}^M\left(\frac{\partial I}{\partial c_i}\right)^2 \right)^{1/2},
\end{eqnarray} which we use MATLAB's \texttt{fminsearch} function to find a root to $\nu$. 
\subsection{Energy-Minimal Controls}
For this, we want to minimize the following functional:
\begin{equation}\label{eq:minenergy}
\min \int_{0}^{t_f}\! u^2 \, dt = \min \int_{x_0}^{x_f}\! u^2 \frac{dt}{dx}\, dx.
\end{equation}
To make \eqref{eq:minenergy} independent of $u$, we note the relation
\begin{equation}\label{eq:dydxtou}
\frac{dy}{dx} = \frac{b_2+a_2y+ux}{b_1+a_1x-uy},
\end{equation} which, solving for $u$, gives
\begin{equation}\label{eq:findu}
u = \frac{-(b_2+a_2y-y'(b_1+a_1x))}{x+yy'}.
\end{equation}
Substituting \eqref{eq:findu} into \eqref{eq:minenergy} gives the integrand 
\begin{equation}\label{eq:energyLagrangian}
L (q,q')= \frac{(b_2+a_2y-y'(b_1+a_1x))^2}{(x+yy')(a_1x^2+b_1x+a_2y^2+b_2y)}.
\end{equation}
This problem can be solved by the Rayleigh-Ritz method as explained in the Time-Minimal section.

\section{Three-dimensional Systems}
Next, instead of considering $N=1$ we allow an arbitrary number of Lindblad operators $N$. We want to find optimal controls for the system \eqref{eq:bloch} with boundary values given by \eqref{eq:endpoint} where $q=[x;y;z]^T$ is in the Bloch ball, that is $x^2+y^2+z^2\leq 1$.
\subsection{Time-Minimal Controls}
This situation is similar to the two-dimensional case. All we need to do is modify \eqref{eq:lagrangian} to
\begin{equation}\label{eq:3dtimelagrangian}
I(q,q') = \int\! L(q,q')  dt =  \int \frac{x+yy'+zz'}{\langle q,b+Bq\rangle}\, dx,
\end{equation}which can be solved with the same algorithm to the two-dimensional case with the following form
\begin{equation}\label{eq:threeguessfunctions}
\begin{split}
&y^i(x) = (x-x_0)(x-x_f)^i,\hbox{ } y^0(x) = \frac{y_f-y_0}{x_f-x_0}(x-x_0)+y_0,\\
&z^i(x) = y^i(x),\quad z^0(x) = \frac{z_f-z_0}{x_f-x_0}(x-x_0)+z_0,
\end{split}
\end{equation} for $i=1,\ldots,M$, where we now have to solve for $2M$ variables.

\subsection{Energy-Minimal Controls}
We want to minimize the cost functional
\begin{equation}\label{eq:3dimenergy}
I = \int_{0}^{t_f} \!\left( \sum_{i=1}^3 \! u_i^2 \right)\, dt = \int_{x_0}^{x_f} \! \left(\sum_{i=1}^3 \! u_i^2 \frac{dt}{dx}\right) \, dx.
\end{equation}
To solve this, we need to make \eqref{eq:3dimenergy} independent of the $u_i$'s. This yields the following system of equations coming from equations \eqref{eq:bloch}:
\begin{equation}\label{eq:systemfor3energy}
\begin{split}
\frac{dy}{dx} & =  \frac{b_2+a_2y-u_3x+u_1z}{b_1+a_1x+u_3y-u_2z}, \\
\frac{dz}{dx} & =  \frac{b_3+a_3z+u_2x-u_1y}{b_1+a_1x+u_3y-u_2z}, \\
\frac{dz}{dy} & =  \frac{b_3+a_3z+u_2x-u_1y}{b_2+a_2y-u_3x+u_1z}.
\end{split}
\end{equation}
\indent We consider $u_1=0$ and we drop the third equation in \eqref{eq:systemfor3energy}. One can alternatively choose $u_2=0$ or $u_3=0$ and the others two controls different to zero. Denoting by $\displaystyle{\Gamma(q,q')=\frac{1}{x^2 + xyy' + xzz'}}$, solving for $u_2$ and $u_3$ yields:

\begin{align}
u_2=&-\Gamma(q,q')\left(b_3x - a_1x^2z' - a_2y^2z'+ a_3xz\right.\nonumber\\
&\left.\qquad\qquad \quad- b_1xz' + b_3yy' - b_2yz' + a_3yy'z\right)\label{eq:u1iszerocontrols},\\
u_3 =&\Gamma(q,q')\left(b_2x - a_1x^2y' - a_3y'z^2 + a_2xy \right.\nonumber\\
&\left.\qquad\quad\qquad-b_1xy' - b_3y'z+b_2zz' + a_2yzz'\right)\label{eq:u1iszerocontrols1}.
\end{align}

To minimize \eqref{eq:3dimenergy}, we must use \eqref{eq:u1iszerocontrols} and \eqref{eq:u1iszerocontrols1} (along with $u_1=0$), the guesses \eqref{eq:threeguessfunctions}, and solve for the coefficients.

\section{Numerical Results}
\subsection{Two-dimensional systems}
We will work an example with parameter values $b_1=1$, $b_2=2$, $a_1=-3$ and $a_2=-4$. Additionally, we will take $\varepsilon = \delta = 10^{-3}$. \\
\indent Solving for the apogee \eqref{eq:endpoint} in polar coordinates, we get $q_{apogee} = [0.4079,0.4493]^T$. Figure $1$ shows a simulation of the trajectory for the $7^{th}$ order curve.

\begin{table}[h!]
\centering
\begin{tabular}{|c||l|l||l|l|}
	\hline
	\multirow{2}{*}{$M$} & \multicolumn{2}{c||}{Time-Minimal} & \multicolumn{2}{|c|}{Energy-Minimal} \\
	\cline{2-5}
	& Time & Energy & Time & Energy \\
	\hline
	$1$ & 1.9371 & 7.5830 & 1.9393 & 0.5365 \\
	\hline
	$3$ & 1.9366 & 8.6873 & 2.1477 & 0.2410 \\
	\hline
	$5$ & 1.9361 & 1.6368 & 2.1789 & 0.2334 \\
	\hline
	$7$ & 1.9359 & 1.3765 & 2.1569 & 0.2369\\
	\hline
\end{tabular}
\caption{Numerical results from time-minimal and energy-minimal controls with solutions of various orders.}
\label{table:2d}
\end{table}

\begin{figure}[h!]
	\includegraphics[scale=0.28]{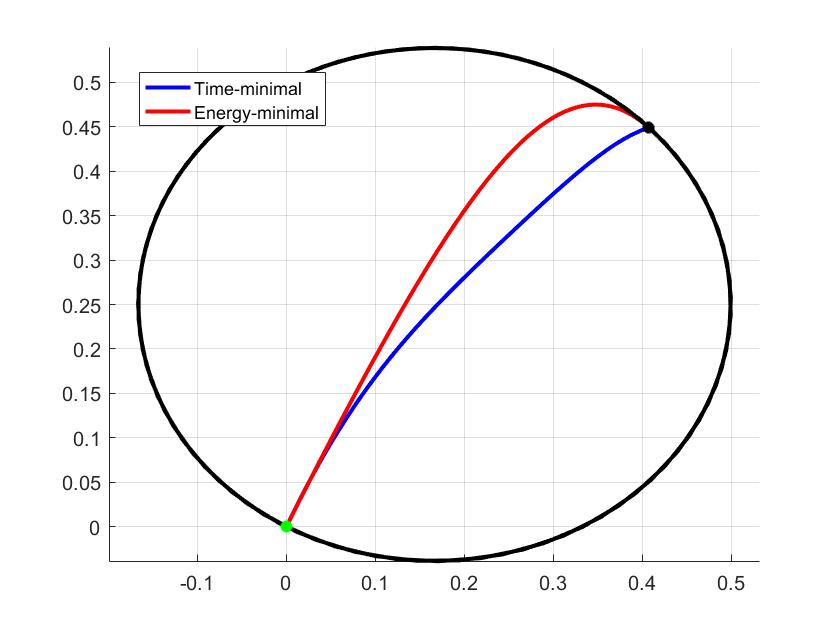}
	\caption{Trajectory of the $7^{th}$ order curve. The black ellipse is the escape chimney.}
\end{figure}

We can also use \eqref{eq:findu} to determine what controls are required for the desired trajectory. Since all the computations are done independent of time, we will report a plot of $u$ versus $x$ in Figure $2$.

\begin{figure}[h!]
	\includegraphics[scale=0.28]{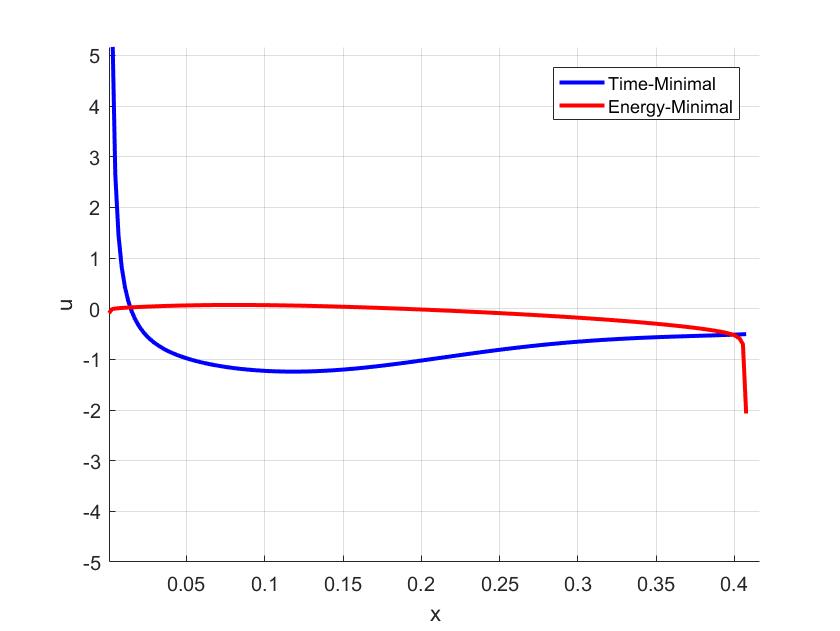}
	\label{figure2}
	\caption{Controls for the $7^{th}$ order curve.}
	\label{fig:2d}
\end{figure}


An expected downfall of the simulations is that multiple local minimal solutions might exist. To find the unique global minimal solution, we repeat the algorithm for various initial conditions. We then find the best solution out of all of the candidates. For this case, we ran the algorithm 25 times and the initial conditions were uniformly randomly chosen in the $l^{\infty}$-ball of radius 2.

\subsection{Three-dimensional systems}
The parameters chosen here will be $b=[1,2,3]^T$ and $B=\text{diag}(-7,-6,-5)$. 
Again, we will take $\varepsilon=\delta=10^{-3}$.\\
\indent Solving for the apogee \eqref{eq:endpoint} in spherical coordinates, we get $q_{apogee} = [0.1140,0.2954,0.6287]^T$. Figure $3$ shows a simulation of the trajectory for the $4^{th}$ order curve. We will follow the same method to avoid local minimums as in the two dimensional case: for all simulations, we solve for the optimal trajectory off of 50 random initial conditions in the $l^{\infty}$-ball of radius 2.

\begin{table}[h!]
	\centering
	\begin{tabular}{|c||l|l||l|l|}
		\hline
		\multirow{2}{*}{$M$} & \multicolumn{2}{c||}{Time-Minimal} & \multicolumn{2}{|c|}{Energy-Minimal} \\
		\cline{2-5}
		& Time & Energy & Time & Energy \\
		\hline
		$1$ & 1.3188 & 207.26 & 1.3243 & 36.365 \\
		\hline
		$2$ & 1.3188 & 47.519 & 1.3205 & 32.491 \\
		\hline
		$3$ & 1.3189 & 42.431 & 1.3212 & 29.356 \\
		\hline
		$4$ & 1.3188 & 49.693 & 1.3214 & 31.682 \\
		\hline
	\end{tabular}
	\caption{Numerical results from time-minimal and energy-minimal controls with solutions of various orders.}
	\label{table:3d}
\end{table}

\begin{figure}[h!]
	\includegraphics[scale=0.45]{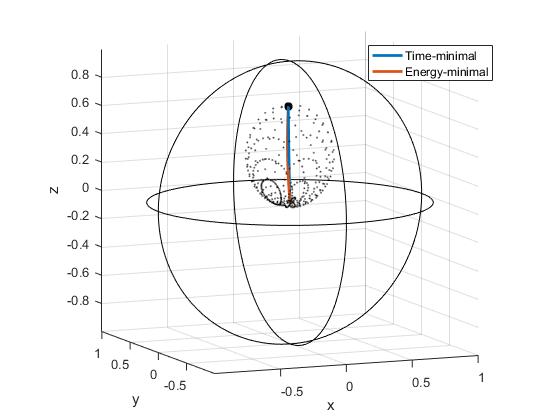}
	\caption{Controls for the $4^{th}$ order curves. The blue curve is the time-minimal curve and the red curve is the energy-minimal. Additionally, the blue ellipsoid is the escape-chimney.}
\end{figure}

\begin{figure}[h!]\label{figureball}
	\includegraphics[scale=0.4]{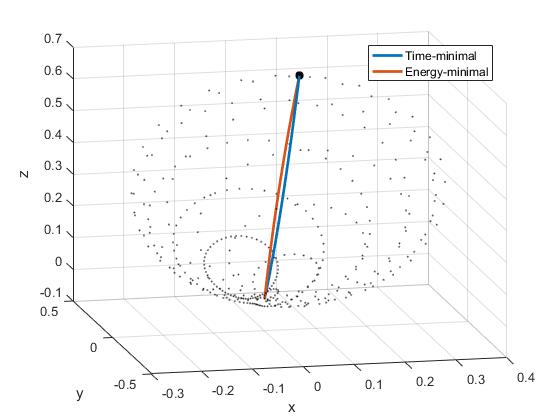}
	\caption{Close up of the escape chimney for the $4^{th}$ order curves.}
\end{figure}

\begin{figure}[h!]
	\includegraphics[scale=0.45]{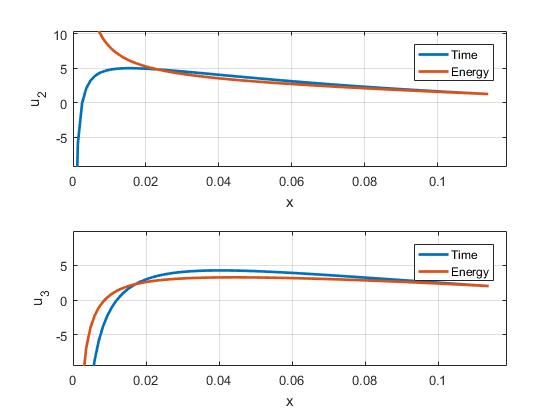}
	\caption{Controls for the $4^{th}$ order curves.}
	\label{fig:3d}
\end{figure}

As before, we can also use \eqref{eq:u1iszerocontrols} and \eqref{eq:u1iszerocontrols1} to determine what controls are required for the desired trajectory. Since all the computations are done independent of time, we will report a plot of $u$ versus $x$ in Figure \ref{fig:3d}.

An interesting feature of the 3-dimensional case is the graphical solutions are much closer than the results for the 2-dimensional case.


Another interesting observation is from Fig. \ref{fig:2d} where the controls approach different values as $x\rightarrow x_f$. In the time-minimal case, $u_{time}(x_f) = -0.5000$. If this value of $u$ is held constant and plugged into \eqref{eq:twodimenions}, then the (stable) fixed point of the system is precisely the apogee of the escape chimney. So under these controls, not only will we approach the apogee but we will also remain there. This also adds intuition to the time-minimal case: the optimal controls make the apogee the stable fixed point.

Now, for the three dimensional case, the controls do approach the same values at the endpoint as seen in Fig. \ref{fig:3d}. These final controls, however, do not make the  apogee the fixed point under \eqref{eq:bloch} i.e. $u_2(x_f)=1.265$ and $u_3(x_f)=2.007$ for both energy and time-minimal trajectories and $-(B+\hat{u})^{-1}b = [0.1364,0.3789,0.5655]^T \ne q_{apogee}$.
This situation will be studied in future work.

\section{Conclusions and future research}
We studied time-minimum and energy-minimum global optimal control problems for dissipative open quantum systems where the dynamics is described by the Lindblad equation and controls are unbounded. We have transformed such a control system into a bi-linear singular control system in the Bloch ball and have come up with the construction of a numerical algorithm to design optimal paths to achieve a desired point given initial states close to the origin of the Bloch ball in both optimal control problems.

All of the results presented are based on having fast control of the Hamiltonian in \eqref{eq:master}, i.e. unbounded controls. It would be interesting to develop both time and energy-minimal control schemes where the control $u$ is bounded (for example, $\lVert u\rVert_{\infty}\leq 1)$. We are currently working on this problem building on the work of \cite{5288564} and \cite{kirk2004optimal}.

Another problem in the bounded control setting is the fact that determining whether the apogee is asymptotically reachable is not clear. We hope to extend the results from \cite{Brockett1975} to determine when the apogee is asymptotically reachable. Extensions of our results to higher-order dimensional systems is another task to work based in this work. Finally, it would also be interesting to determine the best basis of functions for the Rayleigh-Ritz methods as well as the best order of solutions to use.

\addtolength{\textheight}{-12cm}   

\end{document}